\documentclass{amsart}
\usepackage{amssymb}
\usepackage{amsmath}
\def\version{Commuting Jacobi Operators last changed 3 May 2005 Version 2e}
\usepackage{a4}
\usepackage{color}

\begin{document}
\date{\version}
\newtheorem{theorem}{Theorem}[section]
\newtheorem{lemma}[theorem]{Lemma}
\newtheorem{remark}[theorem]{Remark}
\newtheorem{definition}[theorem]{Definition}
\newtheorem{corollary}[theorem]{Corollary}
\newtheorem{example}[theorem]{Example}
\def\qedbox{\hbox{$\rlap{$\sqcap$}\sqcup$}}
\makeatletter
  \renewcommand{\theequation}{%
   \thesection.\alph{equation}}
  \@addtoreset{equation}{section}
 \makeatother
\date{\version}
\def\BB{\mathcal{B}}
\title
{Manifolds with Commuting Jacobi Operators}
\author{M. Brozos-V\'azquez and P. Gilkey}
\begin{address}{MBY:Department of Geometry and Topology, Faculty of Mathematics, University of
Santiago de Compostela, 15782 Santiago de Compostela, Spain}\end{address}
\begin{email}{mbrozos@usc.es}\end{email}
\begin{address}{PG: Mathematics Department, University of Oregon, Eugene, OR 97403, USA}\end{address}
\begin{email}{gilkey@darkwing.uoregon.edu}\end{email}
\begin{abstract} We characterize Riemannian manifolds of constant sectional curvature in terms of commutation properties
of their Jacobi operators.
\end{abstract}
\keywords{Jacobi operator, Osserman manifold, Jacobi-Tsankov manifold, constant sectional curvature.
\newline 2000 {\it Mathematics Subject Classification.} 53C20}
\maketitle
\section{Introduction}

Let $\nabla$ be the Levi-Civita connection of a connected smooth
$m$-dimensional Riemannian manifold $(M,g)$. Let
$\mathcal{R}(x,y):=\nabla_x\nabla_y-\nabla_y\nabla_x-\nabla_{[x,y]}$
be the associated curvature operator. The Jacobi operator
$J(x):y\rightarrow\mathcal{R}(y,x)x$ is intimately related with
the  underlying geometry of the manifold. One says that $(M,g)$ is
{\it Osserman} if the eigenvalues of $J(x)$ are constant on the
sphere bundle $S(M,g)$ of unit tangent vectors. One says that
$(M,g)$ is a {\it local $2$-point homogeneous space} if the
local isometries of $(M,g)$ act transitively on $S(M,g)$; this
necessarily implies that $(M,g)$ is Osserman. Osserman
\cite{refOss} wondered if the converse held; this has been called
the Osserman conjecture by subsequent authors. Chi \cite{Ch88} and
Nikolayevsky \cite{Ni03,Nikp1,Nikp2} established the Osserman
conjecture if $m\ne16$. This gives a very pretty characterization
of local $2$-point homogeneous spaces in terms of the geometry of
the Jacobi operator.

Tsankov \cite{Y05} began a similar characterization of manifolds of constant sectional curvature in terms of the Jacobi operator by showing:
\begin{theorem}\label{thm-1.1}
Let $(M,g)$ be a hypersurface in $\mathbb{R}^{m+1}$. If $J(x)J(y)=J(y)J(x)$ provided $x\perp y$ and if $m\ge3$, then
$(M,g)$ has constant sectional curvature.
\end{theorem}

In this paper, we remove the hypothesis that $(M,g)$ is a hypersurface to characterize flat manifolds and manifolds of constant sectional
curvature in terms of commutation properties of the Jacobi operator:
\begin{theorem}\label{thm-1.2}
Let $(M,g)$ be a Riemannian manifold of dimension $m\ge3$.\begin{enumerate}
\item If $J(x)J(y)=J(y)J(x)$ for all $x,y$, then $(M,g)$ is flat.
\item If $J(x)J(y)=J(y)J(x)$ provided $x\perp y$, then $(M,g)$ has constant sectional curvature.
\end{enumerate}\end{theorem}

It is convenient to work in the purely algebraic context. Let $V$ be a real vector space of dimension $m\ge3$ which is
equipped with a positive definite inner product $\langle\cdot,\cdot\rangle$. Let $R\in\otimes^4V^*$ be an algebraic curvature tensor on $V$,
i.e. $R$ has the usual symmetries of the curvature tensor:
\begin{equation}\label{eqn-1.a}\begin{array}{l}
R(x,y,z,w)=R(z,w,x,y)=-R(y,x,z,w),\\
R(x,y,z,w)+R(y,z,x,w)+R(z,x,y,w)=0\,.\vphantom{\vrule height 11pt}
\end{array}\end{equation}
Let $\mathcal{R}(x,y)$ be the associated {\it curvature operator}; $\mathcal{R}$ is characterized by the identity:
$$\langle\mathcal{R}(x,y)z,w\rangle=R(x,y,z,w)\,.$$
The {\it Jacobi operator} $J(x)$ is then defined by $J(x):y\rightarrow\mathcal{R}(y,x)x$. One says that $R$ is {\it Osserman} if the
eigenvalues of $R$ are constant on the associated unit sphere $S(V):=\{\xi\in V:\langle\xi,\xi\rangle=1\}$. Motivated by Theorem
\ref{thm-1.1},  one says that $R$ is {\it Jacobi-Tsankov} if
$$X\perp Y\quad\text{implies}\quad J(X)J(Y)=J(Y)J(X)\,.$$

The following two tensors will play a crucial role in our investigation. The tensor $R_0$ of constant sectional curvature $+1$ is defined by
$$R_0(x,y)z=\langle y,z\rangle x-\langle x,z\rangle y\,.$$
Let $\Theta$ be a linear endomorphism of $V$. If $\Theta^2=-\text{id}$ and if $\Theta^*=-\Theta$, then $\Theta$ is said to be a {\it
Hermitian almost complex structure} on
$V$. Such a $\Theta$ exists, of course, if and only if $m$ is even. If $\Theta$ is a Hermitian almost complex structure on $V$, set:
$$
R_\Theta(x,y)z:=\langle\Theta y,z\rangle\Theta x-\langle\Theta x,z\rangle\Theta y-2\langle\Theta
x,y\rangle\Theta z\,.
$$

The tensors $R_0$ and $R_\Theta$ play a crucial role in the proof of the Osserman conjecture and also play a crucial role in our present
analysis. We will derive Theorem \ref{thm-1.2} from the following purely algebraic assertion:
\begin{theorem}\label{thm-1.3}Let $\Theta$ be a Hermitian almost complex structure on $V$. \begin{enumerate}
\item If $J(x)J(y)=J(y)J(x)$ for all $x,y$, then $R=0$.
\item $R_0$ and $R_\Theta$ are Jacobi-Tsankov curvature tensors.
\item If $R$ is Jacobi-Tsankov, then either $R=cR_0$ or $R=cR_\Theta$ for some Hermitian almost complex structure $\Theta$.
\end{enumerate}
\end{theorem}

Here is a brief guide to this paper. In Section \ref{sect-2}, we establish Theorem \ref{thm-1.3}. In Section \ref{sect-3},
we use Theorem
\ref{thm-1.3} to prove Theorem \ref{thm-1.2}.

\section{Algebraic results}\label{sect-2}

Throughout this section, let $R$ be an algebraic curvature tensor on a vector space $V$ of dimension $m$  equipped with a
positive definite inner product $\langle\cdot,\cdot\rangle$.

\begin{proof}[Proof of Theorem \ref{thm-1.3} (1)] Suppose that $J(x)J(y)=J(y)J(x)$ for all $x,y\in V$.  The Jacobi operators
$\{J(x)\}|_{x\in V}$ form a commuting family of self-adjoint operators. Such a family can be simultaneously diagonalized - i.e. we can
decompose
$$
V=\oplus_\lambda V_\lambda\quad\text{where}\quad J(x)\xi=\lambda(x)\xi\ \forall\ \xi\in V_\lambda\,.
$$
Choose $\xi\in V$ and decompose
$\xi=\textstyle\sum_\lambda\xi_\lambda$ for $\xi_\lambda\in
V_\lambda$. Let
$$\mathcal{O}:=\{\xi\in V:\xi_\lambda\ne0\ \forall\ \lambda\}\,.$$
 Since $\mathcal{O}=\cap_\lambda\{V-
V_\lambda^\perp\}$, $\mathcal{O}$ is a dense open subset
of $V$. We have
$$0=J(\xi)\xi=\textstyle\sum_\lambda\lambda(\xi)\xi_\lambda\,.$$
Since $\{\xi_\lambda\}$ is a linearly independent set, this
implies $\lambda(\xi)=0$ for all $\lambda$. Fix unit vectors
$\eta_\lambda\in  V_\lambda$. Then $\lambda(\xi)=\langle
J(\xi)\eta_\lambda,\eta_\lambda\rangle$, so $\lambda$ is a
continuous function of $\xi$. As $\lambda(\cdot)$ vanishes on
$\mathcal{O}$, which is a dense open subset of $V$,
$\lambda(\cdot)$ vanishes identically. Thus $J(x)=0$ for all $x\in
V$; it now follows that $R=0$; see, for example, the discussion in
\cite{G02}.
\end{proof}

\begin{proof}[Proof of Theorem \ref{thm-1.3} (2)] Suppose that $\Theta$ defines $R$ where $\Theta$ is a Hermitian almost complex structure on
$V$.  If $x\perp y$, then $\langle\Theta x,\Theta y\rangle=\langle x,y\rangle=0$. Consequently:
\begin{eqnarray*}
&&J(x)\xi=R(\xi,x)x=\langle\Theta x,x\rangle\Theta \xi-\langle\Theta \xi,x\rangle\Theta x-
2\langle\Theta \xi,x\rangle\Theta x=3\langle \xi,\Theta x\rangle\Theta x,\\
&&\{J(x)J(y)-J(y)J(x)\}z=9\langle z,\Theta
y\rangle\langle\Theta y,\Theta x\rangle\Theta x -9\langle
z,\Theta x\rangle\langle\Theta x,\Theta y\rangle\Theta y=0\,.
\end{eqnarray*}
Thus $R_\Theta$ is Tsankov. Similarly if $R=R_0$ is the algebraic curvature tensor of constant sectional curvature
$+1$, then:
\begin{eqnarray*}
&&J(x)\xi=\langle x,x\rangle\xi-\langle\xi,x\rangle x,\\
&&J(x)J(y)z=J(x)\{\langle y,y\rangle z-\langle z,y\rangle y\}\\
&&\qquad\qquad=\langle y,y\rangle\langle x,x\rangle z-\langle y,y\rangle\langle z,x\rangle x
-\langle z,y\rangle\langle x,x\rangle y+\langle z,y\rangle\langle y,x\rangle x\\
&&\qquad\qquad=\langle y,y\rangle\langle x,x\rangle z-\langle y,y\rangle\langle z,x\rangle x
-\langle z,y\rangle\langle x,x\rangle y=J(y)J(x)\,.
\end{eqnarray*}
Assertion (2) of Theorem \ref{thm-1.3} now follows.
\end{proof}

Before establishing Assertion (3) of Theorem \ref{thm-1.3}, we need a number of technical results.
It is convenient to polarize the Jacobi operator to define a bilinear and self-adjoint operator:
$$
J(x,y):z\rightarrow\textstyle\frac12\mathcal{R}(z,x)y+\frac12\mathcal{R}(z,y)x\,.$$
We note that
\begin{equation}\label{eqn-2.a}
\begin{array}{l}
J(x)=J(x,x),\quad J(x,y)y=-\textstyle\frac12J(y)x,\quad\text{and}\\
J(\cos\theta x+\sin\theta y)=\cos^2\theta J(x)+2\cos\theta\sin\theta J(x,y)+\sin^2\theta J(y)\,.\vphantom{\vrule height 11pt}
\end{array}\end{equation}

 Let $r(x):=\operatorname{Rank}\{J(x)\}$. Since $J(x)x=0$, $r(x)\le m-1$ for all $x\in V$.

\begin{lemma}\label{lem-2.1} Let $0\ne R$ be a Jacobi-Tsankov algebraic curvature tensor.
 Suppose that $r(x)=m-1$ for some $x\in S(V)$.
\begin{enumerate}
\item $J(\cdot)$ has maximal rank on an open dense subset of $V$.
\item $R$ has constant sectional curvature $c\ne0$.
\end{enumerate}
\end{lemma}

\begin{proof}
 We clear the previous notation and let $\mathcal{O}:=\{x\in V:r(x)=m-1\}$. As $\mathcal{O}$ is non-empty,
$\mathcal{O}$ is a dense open subset of $V$. Let $x\in\mathcal{O}$ and let $y\in x^\perp$. Then
$$
  J(x)J(y)x=J(y)J(x)x=0\,.
$$
Consequently $g(J(y)x,J(x)z)=0$ for all $z$. As $\operatorname{range}(J(x))=x^\perp$,
we have $g(J(y)x,z)=0$ if $z\perp x$. Thus $R(x,y,y,z)=0$ if $z\perp x$ and $y\perp x$. This relation holds for $x$ belonging to an open
dense subset of $V$. Consequently it holds for all $x\in V$. Thus if $\{e_i\}$ is an orthonormal basis for $V$ and if $\{i,j,k\}$
are distinct indices, we have $R(e_i,e_j,e_j,e_k)=0$.  Thus if $m=3$, the only non-zero curvatures are
$R(e_i,e_j,e_j,e_i)$.

Suppose $m\ge4$. Let $\ell$ be a fourth distinct index. Polarization yields
\begin{equation}\label{eqn-2.b}
R(e_i,e_j,e_\ell,e_k)+R(e_i,e_\ell,e_j,e_k)=0\,.
\end{equation} We use the relations of Equations
(\ref{eqn-1.a}) and (\ref{eqn-2.b}) to see:
\begin{eqnarray*}
0&=&R(e_i,e_j,e_k,e_\ell)+R(e_i,e_k,e_\ell,e_j)+R(e_i,e_\ell,e_j,e_k)\\
&=&R(e_i,e_j,e_k,e_\ell)-R(e_i,e_k,e_j,e_\ell)-R(e_i,e_j,e_\ell,e_k)\\
&=&R(e_i,e_j,e_k,e_\ell)+R(e_i,e_j,e_k,e_\ell)+R(e_i,e_j,e_k,e_\ell)\\
&=&3R(e_i,e_j,e_k,e_\ell)\,.
\end{eqnarray*}
Thus the only non-zero curvatures are $R(e_i,e_j,e_j,e_i)=c_{ij}$. Consider the new basis
$$
  e_\nu(\theta):=\left\{\begin{array}{rll}
\cos\theta e_i+\sin\theta e_j&\text{if}&\nu=i,\\
-\sin\theta e_i+\cos\theta e_j&\text{if}&\nu=j,\\
e_k&\text{if}&\nu\ne i,j\,.\end{array}\right.
$$
We compute
$
0=R(e_i(\theta),e_k,e_k,e_j(\theta))=\cos\theta\sin\theta\{-c_{ik}+c_{jk}\}\,.
$
It now follows that $c_{ik}=c_{jk}$ and consequently $R$ has constant sectional curvature.\end{proof}

 We assume $r(x)<m-1$ henceforth; thus  there exists $y$ so $x\perp y$ and $J(x)y=0$.

\begin{lemma}\label{lem-2.2}
 Let $R$ be a Jacobi-Tsankov algebraic curvature tensor. Assume that $r(x)<m-1$ $\forall x$.
Let $x\in S(V)$. Choose $y\in S(x^\perp)$ so $J(x)y=0$. Then:\begin{enumerate}
\item $J(y)x=0$ and $J(x)J(y)=0$.
\item $0=J(y)^2+J(x)^2-4J(x,y)^2$, $J(x,y)J(x)=J(y)J(x,y)$, and\newline $J(x)J(x,y)=J(x,y)J(y)$.
\item Let $\{x,z_1,z_2\}$ be an orthonormal set. Suppose that $J(x)z_1=\lambda_1z_1$ and $J(x)z_2=\lambda_2z_2$ where $\lambda_1\ne\lambda_2$.
Then
$J(z_1)z_2=0$.
\item We can choose an orthonormal basis for $V$ so that
$$
J(x)=\left(\begin{array}{l}
A\phantom{0}0\phantom{A}0\\
0\phantom{A}0\phantom{A}0\\
0\phantom{A}0\phantom{A}0\end{array}\right),
J(y)=\left(\begin{array}{ll}
0\phantom{A}0\phantom{A}0\\
0\phantom{0}A\phantom{A}0\\
0\phantom{A}0\phantom{A}0\end{array}\right),
J(x,y)=\frac12\left(\begin{array}{l}
0\phantom{0}A\phantom{A}0\\
A\phantom{0}0\phantom{A}0\\
0\phantom{A}0\phantom{A}0\end{array}\right)\,.
$$
\end{enumerate}
\end{lemma}

\begin{proof}
Suppose that $J(x)y=0$. We use the relations of Equation (\ref{eqn-2.a}) to compute:
\begin{equation}\label{eqn-2.c}
\begin{array}{l}
J(\cos\theta x+\sin\theta y)J(-\sin\theta x+\cos\theta y)y\\
\quad=J(\cos\theta x+\sin\theta y)\{\sin^2\theta J(x)-2\sin\theta\cos\theta J(x,y)+\cos^2\theta J(y)\}y\\
\quad=J(\cos\theta x+\sin\theta y)\{\sin\theta\cos\theta J(y)x\}\\
\quad=\{\cos^2\theta J(x)+2\sin\theta\cos\theta J(x,y)+\sin^2\theta J(y)\}\{\sin\theta\cos\theta J(y)x\}\\
\quad=2\sin^2\theta\cos^2\theta J(x,y)J(y)x+\sin^3\theta\cos\theta J(y)J(y)x
\end{array}\end{equation}
and
\begin{equation}\label{eqn-2.d}
\begin{array}{l}
J(-\sin\theta x+\cos\theta y)J(\cos\theta x+\sin\theta y)y\\
\quad=J(-\sin\theta x+\cos\theta y)\{\cos^2\theta J(x)+2\cos\theta\sin\theta J(x,y)+\sin^2\theta  J(y)\}y\\
\quad=\{\sin^2\theta J(x)-2\cos\theta\sin\theta J(x,y)+\cos^2\theta J(y)\}\{-\cos\theta\sin\theta J(y)x\}\\
\quad=2\cos^2\theta\sin^2\theta J(x,y)J(y)x-\cos^3\theta\sin\theta J(y)J(y)x\,.
\end{array}\end{equation}
Subtracting Equation (\ref{eqn-2.d}) from Equation (\ref{eqn-2.c})
yields $\sin\theta\cos\theta J(y)^2x=0$ and hence as $\theta$ was
arbitrary, $J(y)^2x=0$. Since $J(y)$ is diagonalizable,
$J(y)x=0$.

We have $J(x)J(y)y=0$ and $J(x)J(y)x=0$. Let $z\perp\{x,y\}$. To complete the proof of Assertion (1), we must show $J(x)J(y)z=0$. We compute:
\begin{eqnarray*}
0&=&J(\cos\theta x+\sin\theta z)J(y)x\\
&=&J(y)J(\cos\theta x+\sin\theta z)x\\
&=&J(y)\{\cos^2\theta J(x)+2\cos\theta\sin\theta J(x,z)+\sin^2\theta J(z)\}x\\
&=&2\cos\theta\sin\theta J(y)J(x,z)x+\sin^2\theta J(y)J(z)x\\
&=&-\cos\theta\sin\theta J(y)J(x)z+\sin^2\theta J(z)J(y)x\\
&=&-\cos\theta\sin\theta J(y)J(x)z\,.
\end{eqnarray*}

We now prove Assertion (2). Because $J(x,y)x=-\frac12J(x)y=0$ and because $J(x,y)y=-\frac12J(y)x=0$, we have:
\begin{eqnarray*}
&&J(x,y)\{-\sin\theta x+\cos\theta y\}=0,\quad\text{so}\\
&&J(\cos\theta x+\sin\theta y)\{-\sin\theta x+\cos\theta y\}=0\,.
\end{eqnarray*}
Thus by applying Assertion (1) to the pair $\{\cos\theta x+\sin\theta y,-\sin\theta x+\cos\theta y\}$,
\begin{eqnarray*}
0&=&J(\cos\theta x+\sin\theta y)J(-\sin\theta x+\cos\theta y)\\
&=&\{\cos^2\theta J(x)+2\sin\theta\cos\theta J(x,y)+\sin^2\theta J(y)\}\\
  &&\cdot \{\sin^2\theta J(x)-2\sin\theta\cos\theta J(x,y)+\cos^2\theta J(y)\}\\
&=&\cos^2\theta\sin^2\theta\{J(y)^2+J(x)^2-4J(x,y)^2\}\\
&+&2\sin^3\theta\cos\theta\{J(x,y)J(x)-J(y)J(x,y)\}\\
&+&2\sin\theta\cos^3\theta\{J(x,y)J(y)-J(x)J(x,y)\}\,.
\end{eqnarray*}
Assertion (2) now follows since this identity holds for all $\theta$.

Let $\{x,z_1,z_2\}$ be an orthonormal set with $J(x)z_i=\lambda_iz_i$ where $\lambda_1\ne\lambda_2$. To prove Assertion (3), we compute
\begin{eqnarray*}
&&J(x)J(\cos\theta z_1+\sin\theta z_2)z_1=J(x)\{2\cos\theta\sin\theta J(z_1,z_2)+\sin^2\theta J(z_2)\}z_1\\
&&\qquad=J(x)\{-\cos\theta\sin\theta J(z_1)z_2+\sin^2\theta J(z_2)z_1\}\\
&&\qquad=-\lambda_2\cos\theta\sin\theta J(z_1)z_2+\lambda_1\sin^2\theta J(z_2)z_1,\\
&&J(\cos\theta z_1+\sin\theta z_2)J(x)z_1=\lambda_1\{2\cos\theta\sin\theta J(z_1,z_2)+\sin^2\theta J(z_2)\}z_1\\
&&\qquad=-\lambda_1\cos\theta\sin\theta J(z_1)z_2+\lambda_1\sin^2\theta J(z_2)z_1\,.
\end{eqnarray*}
Assertion (3) now follows since we have
$$\lambda_2J(z_1)z_2=\lambda_1J(z_1)z_2\,.$$

To prove the final assertion,  choose a basis $\{e_1,...,e_r\}$ for $\operatorname{Range}(J(x))$ so that
$$J(x)e_i=\lambda_ie_i\quad\text{for}\quad\lambda_i\ne0\,.$$
We then have $J(y)e_i=0$ and thus $4J(x,y)^2e_i=\lambda_i^2e_i$. We set
$$
f_i:=2\lambda_i^{-1}J(x,y)e_i\,.
$$
 We then compute:
\begin{eqnarray*}
&&\langle f_i,f_j\rangle=4\lambda_i^{-1}\lambda_j^{-1}\langle J(x,y)e_i,J(x,y)e_j\rangle=
4\lambda_i^{-1}\lambda_j^{-1}\langle
J(x,y)^2e_i,e_j\rangle =\delta_{ij},\\
&&J(x)f_i=2\lambda_i^{-1}J(x)J(x,y)e_i=2\lambda_i^{-1}J(x,y)J(y)e_i=0\quad\text{so}\\
&&f_i\in\ker(J(x))=\operatorname{Range}(J(x))^\perp=\operatorname{Span}\{e_1,...,e_r\}^\perp\,.
\end{eqnarray*}
Thus $\{e_1,...,e_r,f_1,...,f_r\}$ is an orthonormal set. Set
$$A=\operatorname{diag}(\lambda_1,....,\lambda_\ell)\,.$$
Since $J(y)J(x,y)=J(x,y)J(x)$, we have $J(y)f_i=\lambda_if_i$. Since
$J(x,y)e_i=\frac12\lambda_if_i$ and $J(x,y)f_i=2\lambda_i^{-1}J(x,y)^2e_i=\frac12\lambda_ie_i$, on
$V_i:=\operatorname{Span}\{e_i,f_i\}$ we then have on $\operatorname{Span}\{e_1,...,e_\ell,f_1,...,f_\ell\}$ that
$$
J(x)=\left(\begin{array}{lll}
A&0\\
0&0\end{array}\right),\quad
J(y)=\left(\begin{array}{lll}
0&0\\
0&A\end{array}\right),\quad
J(x,y)=\left(\begin{array}{lll}
0&\frac12A\\
\frac12A&0\end{array}\right)\,.
$$

{ Let
$W=\operatorname{Span}\{e_1,...,e_\ell,f_1,...,f_\ell\}^\perp$.
Since $\operatorname{Span}\{e_1,...,e_\ell,f_1,...,f_\ell\}$ is
preserved by $\{J(x),J(y),J(x,y)\}$, $W$ is preserved by these
operators. As $W\subset\operatorname{Range}\{J(x)\}^\perp$,
$J(x)=0$ on $W$. Thus $J(y)^2=4J(x,y)^2$ on $W$. If $J(y)z=\mu z$
for $z\in W$, then
$$
\mu^3z=J(y)^3z=\textstyle\frac14J(x,y)^2J(y)z=\textstyle\frac14J(x,y)J(x)J(x,y)z
\,.$$
Since $J(x,y)z\in W$, $J(x)J(x,y)z=0$ and thus $\mu=0$. This shows $J(y)$ and hence $J(x,y)$ vanishes on $W$. Choosing
an orthonormal basis for $W$ then leads to the desired decomposition.}
\end{proof}

We continue our study. Let
\begin{eqnarray*}
&&W(x):=\mathbb{R}\cdot x\oplus\operatorname{Range}(J(x))\,.
\end{eqnarray*}
\begin{lemma}\label{lem-2.3}
 Let $0\ne R$ be a Jacobi-Tsankov algebraic curvature tensor. Assume  that $r(x)<m-1$ for all $x\in V$.
Let $x\in S(V)$ with $r(x)\ne0$. If $w\in S(W(x))$,
\begin{enumerate}
\item $\operatorname{Range}(J(w))\subset W(x)$ and $J(w)$ vanishes on $W(x)^\perp$.
\item $J(w)$ is similar to $J(x)$.
\item $J(x)$ has exactly two eigenvalues.
\end{enumerate}\end{lemma}

\begin{proof} Fix $w\in S(W(x))$. Expand $w=a_0x+\textstyle\sum a_iw_i$ where
$J(x)w_i=\lambda_iw_i$ for
$\lambda_i\ne0$. Let $y\in S(W(x)^\perp)$. We apply Lemma \ref{lem-2.2}. As $y\perp\operatorname{Range}(J(x))$,
$J(x)y=0$ so $J(y)x=0$. Furthermore since $J(x)y=0$, since $J(x)w_i=\lambda_iw_i$, and since $\lambda_i\ne0$, $J(y)w_i=0$.
Thus $J(y)w=0$ and consequently $J(w)y=0$ for all
$y\in W(x)^\perp$. Thus
$$\operatorname{Range}(J(w))\subset W(x)\quad\text{and}\quad J(w)=0\quad\text{on}\quad W(x)^\perp\,.$$
This proves Assertion (1). Furthermore
$J(x)y=0$ and $J(w)y=0$ implies $J(x)$ is similar to $J(y)$ and $J(w)$ is similar to $J(y)$. This establishes Assertion (2).

To show that Assertion (3) is true, we apply Assertion (2) to see
$$\operatorname{Rank}(J(w))=\operatorname{Rank}(J(w)|_{W(x)})=\dim(W(x))-1=r(x)\,.$$
Suppose $J(x)$ has two distinct non-zero eigenvalues
$\lambda_i\ne\lambda_j$. Then $J(w_i)w_j=0$. Since $J(w_i)w_i=0$,
we would have $\operatorname{Rank}\{J(w_i)\}\le r(x)-1$ which is
false. Thus $J(x)=\lambda\operatorname{id}$ on
$\operatorname{Range}(J(x))$. This shows that $J(x)$ has at most
$2$ eigenvalues.

To complete the proof, we must show $\lambda\ne0$. Suppose to the
contrary that $\lambda=0$. This means $J(x)=0$. Let $z\in S(V)$.
Choose $y\in S(x^\perp\cap z^\perp)$. Then $J(x)y=0$ so $J(x)$ is
similar to $J(y)$ and thus $J(y)=0$. Since $J(y)z=0$, $J(y)$ is
similar to $J(z)$ and $J(z)=0$. This implies $J\equiv0$ and hence
$R=0$ which is false. Thus $J(x)\ne0$ and $\lambda\ne0$.
\end{proof}

\begin{lemma}\label{lem-2.4}
 Let $0\ne R$ be a Jacobi-Tsankov algebraic curvature tensor. Assume that $r(x)<m-1$ for all $x\in V$.
\begin{enumerate}\item $r(x)=1$.
\item $R$ is Osserman.
\end{enumerate}
\end{lemma}

\begin{proof}  Choose $y\in S(x^\perp)$ with $J(x)y=0$. Let
$e_0=x$ and
$f_0=y$. Let
$\lambda$ be the non-zero eigenvalue for
$J(x)$ and let
$r=r(x)$. Choose an orthonormal basis $\{e_0,...,e_r,f_0,...,f_r,g_1,...,g_\ell\}$
for $V$ so that
$$\begin{array}{lll}
J(e_0)e_j=\lambda(1-\delta_{0j})e_j&J(e_0)f_j=0,&J(e_0)g_k=0,\\
J(f_0)f_j=\lambda(1-\delta_{0j})f_j&J(f_0)e_j=0,&J(f_0)g_k=0,\vphantom{\vrule height 11pt}\\
J(e_0,f_0)e_j=\textstyle\frac12\lambda(1-\delta_{0j})f_j&
J(e_0,f_0)f_j=\textstyle\frac12\lambda(1-\delta_{0j})e_j&
J(e_0,f_0)g_k=0\,.\vphantom{\vrule height 11pt}
\end{array}$$
 As $J(e_1)$ preserves $\operatorname{Span}\{e_0,...,e_r\}=W(e_0)$, as $\lambda$ is an eigenvalue of multiplicity $r$
for $J(e_1)$ on $W(e_0)$, that $J(e_1)$ vanishes on $W(e_0)^\perp$, and as $J(e_1)e_1=0$,
$$
J(e_1)e_j=\lambda(1-\delta_{1j})e_j,\quad J(e_1)f_j=0,\quad J(e_1)g_k=0\,.
$$
Let $\xi=\frac1{\sqrt{2}}(e_0+f_0)$. Then $J(\xi)=\frac12\{J(e_0)+J(f_0)+2J(e_0,f_0)\}$.  We show $r=1$ and prove
Assertion (1) by deriving the following contradiction:
\begin{eqnarray*}
&&J(\xi)e_2=\textstyle\frac12\lambda(e_2+f_2),\quad J(\xi)f_2=\textstyle\frac12\lambda(e_2+f_2),\\
&&J(\xi)J(e_1)f_2=0,\quad\qquad J(e_1)J(\xi)f_2=\textstyle\frac12\lambda^2e_2\,.
\end{eqnarray*}

Fix $e\in S(V)$. Consider the $2$-plane
$$\pi:=\operatorname{Span}\{e,\operatorname{Range}(J(e))\}\,.$$
Decompose $x\in S(V)$ in the form $x=\cos\theta e_1+\sin\theta f_1$ for $\theta\in[0,\frac\pi2]$, $e_1\in S(\pi)$ and
$f_1\in S(\pi^\perp)$; $e_1$ is not unique if $\theta=\frac\pi2$ and $f_1$ is not unique if $\theta=0$.
As
$e_1\in\pi$,
$\operatorname{Range}(J(e_1))\subset\pi$ so
$$\pi=\operatorname{Span}\{e_1,\operatorname{Range}(J(e_1))\}\,.$$
Since $f_1\perp\pi$, Lemma \ref{lem-2.2} pertains. As
$J(f_1)$ is similar to $J(e_1)$ and also to $J(e)$, $\lambda(e_1)=\lambda(f_1)=\lambda(e)$. By Lemma \ref{lem-2.2}, we
can extend $\{e_1,f_1\}$ to an orthonormal set
$\{e_1,e_2,f_1,f_2\}$ so that
$$\begin{array}{lll}
J(e_1)e_2=\lambda(e) e_2,&J(f_1)e_2=0,&J(e_1,f_1)e_2=\textstyle\frac12\lambda(e) f_2,\\
J(e_1)f_2=0,& J(f_1)f_2=\lambda(e)
f_2,&J(e_1,f_1)f_2=\textstyle\frac12\lambda(e)
e_2\,.\vphantom{\vrule height 11pt}
\end{array}$$
We may now compute
\begin{eqnarray*}
&&J(x)(\cos\theta e_2+\sin\theta f_2)\\
&=&(\cos^2\theta J(e_1)+2\cos\theta\sin\theta J(e_1,f_1)+\sin^2\theta J(f_1))(\cos\theta e_2+\sin\theta f_2)\\
&=&\lambda(e)(\cos^3\theta e_2+\cos^2\theta\sin\theta f_2+\cos\theta\sin^2\theta e_2+\sin^3\theta f_2)\\
&=&\lambda(e)(\cos\theta e_2+\sin\theta f_2)\,.
\end{eqnarray*}
This implies that $\lambda(x)=\lambda(e)$ as desired.
\end{proof}

\begin{proof}[Proof of Theorem \ref{thm-1.3} (3)] Suppose $0\ne R$ is Jacobi-Tsankov. If $r(x)=m-1$ for any vector $x\in S(V)$, then $R$ has
constant sectional curvature by Lemma \ref{lem-2.1}. On the other hand, if $r(x)<m-1$ for every point $x\in S(V)$, then $R$ is Osserman and
$\operatorname{Rank}\{J(x)\}=1$ for every $x\in S(V)$ by Lemma \ref{lem-2.4}. Work of Chi \cite{Ch88} then shows $R=R_\Theta$ for some
Hermitian almost complex structure on $V$.
\end{proof}

\section{The Geometric Setting}\label{sect-3}
\begin{proof}[Proof of Theorem \ref{thm-1.2}] Let $(M,g)$ be a Riemannian manifold with $m\ge3$. Suppose first that
$J(x)J(y)=J(y)J(x)$ for $x,y\in T_PM$. Then $R_P=0$. Consequently $(M,g)$ is flat; this proves Assertion (1).

Suppose $(M,g)$ is Jacobi-Tsankov and $m\ge3$. If $m$ is not even, then $R_P$ has constant sectional curvature at each point $P$ of $M$ and
hence $(M,g)$ has constant sectional curvature globally. Let $\mathcal{O}$ be open subset of points $P\in M$ so that there exists
a  unit tangent vector $x(P)$ with $r(x(P))=m-1$. Then $g|_{\mathcal{O}}$ has constant sectional curvature. Thus $R=cR_0$ on $\mathcal{O}$
and hence $R=cR_0$ on the closure of $\mathcal{O}$.  Thus $\mathcal{O}$ is an open and closed subset of $M$ and hence all of $M$.

Thus if $(M,g)$ does not have constant sectional curvature,
$r(x)<m-1$ for all $x\in S(M,g)$. Suppose $(M,g)$ is not flat;
there is then a point $P\in M$ and a tangent vector $x\in S(T_PM)$
so that $r(x)=1$. Let $\mathcal{U}$ be a small open
contractable of $P$. Then $r(x)=1$ for $x\in S(\mathcal{U})$. Thus
there is an almost complex structure $\Theta(Q)$ defined on $T_Q$
for every $Q\in U$ so that $R=\lambda R_\Theta$. Furthermore, the
almost complex structure is uniquely determined up to sign. After
a bit of technical fuss, one can see that $\Theta$ can be chosen
to vary smoothly with $Q$, at least locally; global questions are
irrelevant to our argument. The metric in question is Einstein and
thus $\rho(x,x)=\lambda(x)$ is constant. Thus $(M,g)$ is globally
Osserman. The work of Chi \cite{Ch88} then implies $(M,g)$ is a
complex space form - i.e. isometric to $\mathbb{CP}^j$ or its
negative curvature dual. This is false as these manifolds have
$\operatorname{Rank}(J(x))\ge2$. This completes the
proof.\end{proof}

\section*{Acknowledgments}
The research of M. Brozos-V\'azquez was partially supported by project BFM 2003-02949 (Spain).  The research of P. Gilkey was
partially supported by the Max Planck Institute for the Mathematical Sciences (Leipzig, Germany).  It
is a pleasure for both authors to acknowledge helpful conversations with C. Dunn.


\begin{thebibliography}{AAA}

\bibitem{Ch88} Q.-S. Chi,
{\it A curvature characterization of certain locally rank-one symmetric spaces},
J. Differential Geom. {\bf 28} (1988), 187--202.

\bibitem{G02}P. Gilkey, {\bf Geometric Properties of
Natural Operators Defined by the Riemann Curvature Tensor}, World Scientific (2002).

\bibitem{GIZ02} P. Gilkey, R. Ivanova, and T. Zhang, Higher order Jordan Osserman, Pseudo-Riemannian manifolds,
Class and Quantum Gravity {\bf 19} (2002), 4543--4551.

\bibitem{Ni03} Y. Nikolayevsky,
{\it Two theorems on Osserman manifolds},
Differential Geom. Appl. {\bf 18} (2003), 239--253.

\bibitem{Nikp1} Y. Nikolayevsky,
{\it Osserman manifolds of dimension $8$}, {\it Manuscr. Math.} {\bf 115} (2004), 31--53.

\bibitem{Nikp2} Y. Nikolayevsky, Osserman Conjecture in dimension $n \ne 8, 16$; math.DG/0204258.

\bibitem{refOss} R. Osserman,
  {\it Curvature in the eighties},
  Amer. Math. Monthly, {\bf97}, (1990) 731--756.


\bibitem{Y05} Y. Tsankov, A characterization of $n$-dimensional hypersurface in $\mathbb{R}^{n+1}$ with commuting curvature
operators, preprint.

\end{thebibliography}
\end{document}